\documentclass[12pt]{amsart}

\usepackage{amsmath}
\usepackage{amssymb}
\usepackage{latexsym}
\usepackage{amscd}
\usepackage{graphics}
\usepackage{wrapfig}

\renewcommand{\epsilon}{\varepsilon}

\newcommand{\bd}{\partial}

\newcommand{\Z}{\mathbb{Z}}

\newcommand{\img}{\mathrm{Im}\ }
\newcommand{\Ker}{\mathrm{Ker}\ }

\newcommand{\R}{\mathbb{R}}
\newcommand{\Eq}{_{\mathrm{eq}}}

\newcommand{\Til}{\widetilde}
\newcommand{\lm}{\lambda}
\newcommand{\TN}{\Til{N}}
\newcommand{\Tno}{\Til{N_0}}
\newcommand{\Sg}{\Sigma}
\newcommand{\Emb}{\mathrm{Emb}}

\newcommand{\con}{\mathrm{con}}

\newcommand{\ts}{{\otimes 2}}
\newcommand{\e}{\varepsilon}

\newcommand{\picHVD}{
\begin{picture}(12,10)
\put(1,-1){\line(1,0){10}}
\put(1,-1){\line(0,1){10}}
\put(1,9){\line(1,0){10}}
\put(11,-1){\line(0,1){10}}
\put(1,-1){\line(1,1){10}}
\put(1,4){\line(1,0){10}}
\put(6,-1){\line(0,1){10}}
\end{picture}
}
\newcommand{\picVD}{
\begin{picture}(12,10)
\put(1,-1){\line(1,0){10}}
\put(1,-1){\line(0,1){10}}
\put(1,9){\line(1,0){10}}
\put(11,-1){\line(0,1){10}}
\put(1,-1){\line(1,1){10}}
\put(6,-1){\line(0,1){10}}
\end{picture}
}
\newcommand{\picH}{
\begin{picture}(12,10)
\put(1,-1){\line(1,0){10}}
\put(1,-1){\line(0,1){10}}
\put(1,9){\line(1,0){10}}
\put(11,-1){\line(0,1){10}}
\put(1,4){\line(1,0){10}}
\end{picture}
}

\newcommand{\picV}{
\begin{picture}(12,10)
\put(1,-1){\line(1,0){10}}
\put(1,-1){\line(0,1){10}}
\put(1,9){\line(1,0){10}}
\put(11,-1){\line(0,1){10}}
\put(6,-1){\line(0,1){10}}
\end{picture}
}

\begin{document}
\title{
Embedding punctured $n$-manifolds
in Euclidean $(2n-1)$-space}

\author{Dmitry Tonkonog
}

\address{Department of Differential Geometry and Applications, Faculty of Mechanics and Mathematics,
Moscow State University, Moscow, 199992, Russia.}
\email{dtonkonog@gmail.com}

\date{}

\begin{abstract}

Let $N$ be a closed orientable connected
$n$-manifold, 
$n\ge 4$.
We classify
embeddings
of the punctured manifold
$N_0$ into $\R^{2n-1}$
up to isotopy.
Our result in some sense extends results of
J.C.~Becker -- H.H. Glover (1971) and O.~Saeki (1999).

\end{abstract}

\maketitle

\section{Introduction and main results}

This paper is on
the classical Knotting Problem: 
for a manifold $N$ and
a number $m$
describe the set $\Emb^m(N)$ of
isotopy classes of embeddings $N \to \R^m$.
For recent surveys, see
[Sk06, HCEC].

We classify
embeddings
of the punctured $n$-manifold
$N_0$ into $\R^{2n-1}$
up to isotopy.

Unless otherwise stated,
we work in the PL (piecewise linear)
or DIFF (smooth) category
and the results are valid in both categories.
For a manifold $X$ we denote by $\Emb^s X$
the set of isotopy classes of embeddings $X\to\R^s$,
and by $X_0$ we denote $X$ minus the interior of a codimension~0
open ball.

Let $N$ be a closed orientable connected
$n$-manifold.
For $n\ge 6$ the set $\Emb^{2n-1}N$
was described in [Ya84], see also [Sk10].
For $n=3$ and $H_1(N)$ torsion-free,
$\Emb^{2n-1}N_0=\Emb^5N_0$ was described
by Saeki [Sa99], see below. 
The method of [Sa99]
used classification of normal bundles 
of embeddings $N_0\to \R^5$ and could not
be directly generalized on higher dimensions.
Our main result is Theorem~1 below which 
implies a description of
$\Emb^{2n-1}N_0$ for $n\ge 4$, see Corollary~1.

If (co)homology coefficients
are omitted, they are assumed to be $\Z$.
We denote $\Z_{(n)}:=\Z$ when $n$ is even and $\Z_{(n)}:=\Z_2$
when $n$ is odd.
For an abelian group $G$ denote by $B_s(G)$ the group of 
symmetric (for $s$ even) or antisymmetric (for $s$ odd) elements of $G^\ts:=G\otimes_\Z G$.

\smallskip
{\bf Corollary 1.}
{\it 
Let $N$ be a closed connected orientable 
$n$-manifold, $n\ge 4$.
Then, as sets,
$$
\Emb^{2n-1}N_0=A\times H_1(N;\Z_{(n-1)})
$$
where $A$ is a quotient set of $H_1(N)^\ts$.
}

\smallskip
We also present a sketch of proof of the following conjecture,
a detailed proof of which will appear in a subsequent version
of the paper.

\smallskip
{\bf Conjecture 1.}
{\it 
Let $N$ be a closed connected orientable 
$n$-manifold, $n\ge 4$.
Then, as sets,
$$
\Emb^{2n-1}N_0=B_{n}(H_1(N))\times H_1(N;\Z_{(n-1)}).
$$
}

\begin{figure}[h]
\centering
\includegraphics{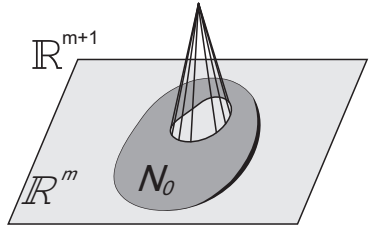}
\caption{The cone map $\Lambda$ which adds a cone to an embedding of $N_0$.}
\label{fig:cone}
\end{figure}

Define the {\it cone map}
$$\Lambda: \Emb^{m}(N_0)\to \Emb^{m+1}(N)$$
which adds a cone over $\bd N_0$, see figure~\ref{fig:cone}.
This map is well-defined
in the PL category and for $m\ge 3n/2+1$ in the smooth category.
\footnote{The sphere $\bd N_0$ is unknotted in $\R^m$ and
$\R^{m+1}$ for $m\ge 3n/2+1$, so we can smoothen the cone by changing a
neighborhood of the cone's vertex.}

\smallskip
{\bf Theorem 1.}
{\it 
Let $N$ be a closed
homologically $k$-connected 
orientable 
$n$-manifold,
$n\ge 4k+4$, $k\ge 0$.
Then the cone map
$\Lambda:\Emb^{2n-2k-1}N_0\to\Emb^{2n-2k}N$
is surjective and the preimage $\lambda^{-1}x$ of each element
is in 1-1 correspondence with
a quotient set of
$H_{k+1}(N)^\ts$. 

In other words, there is the following exact sequence of sets with an action $a$.
$$
\begin{CD}
H_{k+1}(N)^\ts @>a>>\Emb^{2n-2k-1}N_0@>\Lambda>>\Emb^{2n-2k}N@>>>0.
\end{CD}
$$

}

The surjectivity of $\Lambda$ in Main Theorem~1 was known,
see Theorem~2 below
(compare [Vr89, Corollary 3.3]).
Our main result in Theorem~1 is the estimation
of the `kernel' of the cone map $\Lambda$.

\smallskip
{\bf The Becker--Glover Theorem 2.} [BG71]
{\it
Let $N$ be a closed homologically $k$-connected $n$-manifold
and $m\ge 3n/2+2$.
The cone map $\Lambda: \Emb^m (N_0)\to\Emb^{m+1}(N)$
is one-to-one
for
$m\ge 2n-2k$
and is surjective
for $m=2n-2k-1$.
}

\smallskip
Analogously to 
Theorem 1
we could prove that
there is an exact sequence with an action $a$
$$
\begin{CD}
B_3(H_{1}(N))@>a>>\Emb^{5}N_0@>>>H_1(N,\Z)@>>>0.
\end{CD}
$$
in the DIFF category.
This result for
$H_1(N)$ torsion free is known [Sa99].
Moreover, in this case each stabilizer of the action $a$ is in 1-1 correspondence with 
$I_N:=\{\iota_a(\mu):=\cdot\cap\cdot\cap a\ |\ a\in H_1(M)\}\subset B_3(H_1(N))$
where $\mu:H_1(N)\times H_1(N)\times H_1(N)\to \Z$
is the intersection form.
In other words, the preimage of $\Lambda$
is in 1-1 correspondence with $H_1(N)\times B_3(H_1(N))/I_N$.



\smallskip
We conjecture there is a geometric construction
of the invariant $\Emb^{2n-1}N_0\to H_1(N)^\ts$
implicitly obtained in Corollary~1.
If $N$ is a connected sum of several $(S_1\times S^{n-1})$'s
then the 
pairwise linking coefficients $\mathrm{lk} (*\times S^{n-1}_i, *\times S^{n-1}_j)$
provide such an invariant if we take into account
that $H_1(N)\cong H_{n-1}(N)$ in this case 
(the invariant will actually be in $B_n(H_1(N))$).
In the general case, even if $H_1(N)$ is free, 
the linking coefficient is harder to define since $(n-1)$-
cycles of $N_0$ may intersect each other.

\smallskip
{\it Proof of Corollary~1}.
Corollary~1 follows from Theorem~1(a)
and the fact there
is a bijection 
$\Emb^{2n}N\to H_1(N,\Z_{(n-1)})$,
the {\it Whitney invariant} [EBSR].
$\blacksquare$

\section{Proof of Theorem~1}

{\it Proof of Theorem 1 for $k=0$.}
Denote $m=2n-2k$. Consider the following commutative diagram
of sets, in which the horizontal maps are bijections.
$$
\begin{CD}
\Emb^{m}(N)
@.
\begin{picture}(1,10)
\put(-65,3){\vector(1,0){130}}
\end{picture}
^{\alpha}
@.
\pi\Eq^{m-1}(\Til{N})\\
@A\Lambda AA
@.
@A\lambda^\star AA\\
\Emb^{m-1}(N_0)@>\alpha_0>>
\pi\Eq^{m-2}(\Til{N_0})@>\Sigma>>\pi\Eq^{m-1}(\Sigma \Til{N_0})\end{CD}
$$
Here

$\bullet$ $\Til N$ is the {\it deleted product} of $N$,
i.e. $N^2$ minus an open tubular neighborhood
of the diagonal, with standard involution; 

$\bullet$ $\pi^i\Eq X$ is the set of equivariant maps $X\to S^i$
up to equivariant homotopy.

$\bullet$ maps $\alpha$ and $\alpha_0$ are the {\it Haefliger-Wu invariants} [Sk06, 5.2].

$\bullet$ $\Sigma$ is the suspension; 

$\bullet$ $\lambda^*$ is induced by
an equivariant map $\lambda: \Til{N}\to\Sigma \Til{N_0}$ defined below.


\smallskip
{\it Construction of $\lambda$.}
\begin{figure}[h]
\centering
\includegraphics{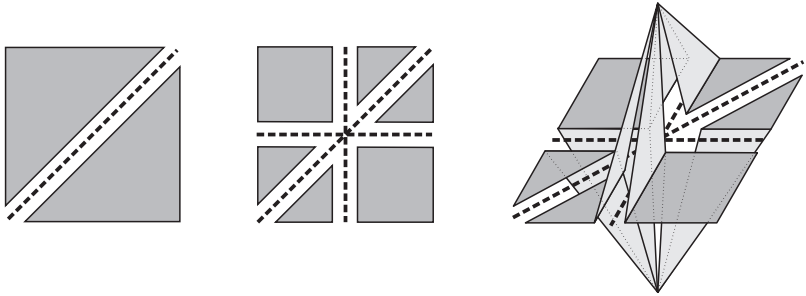}
\caption{From left to right: $\Til {N}$, $\Til {N_0}$ and the image $\lambda\Til{N}\subset\Sigma\Til{N_0}$.}
\label{fig:X1N0}
\end{figure}
We repeat the construction of [BG71].
Represent
\begin{equation}
\label{eqn1}
\Sigma\Til{N_0}=\frac{\Til{N_0}\times[-1;1]}{\Til{N_0}\times\{-1\},\ \Til{N_0}\times\{1\}}.
\end{equation}
For $x\in \Til{N_0}$ set $\lambda(x):=(x,0)$.
We identify $U_\e(P)$ with the unit ball in $\R^n$,
with $P$ corresponding to $0\in \R^n$.
Now set (see figure~2)
$$\lambda(x):=((x_1,v),t-1)
\quad\mathrm{for}\quad x=(x_1,tv)\in N_0\times U_\e(P)
\quad\mbox{where} \quad
x_1\in N_0;
\quad
v\in \bd U_\e(P);
\quad
t\in [0;1].$$
Analogously, for
$x=(tv,x_1)\in B^n\times N_0$ set
$\lambda(x):=((x_1,v),1-t)$.


\smallskip
{\it Proof of commutativity of the diagram above.}
Consider an embedding $f: N_0\subset\R^{m-1}$.
It induces an equivariant map $f_*: \Til{N_0}\to S^{m-2}$.
By definition of the Haefliger-Wu invariant, $[f_*]=\alpha_0[f]$.
\footnote{
Square brackets denote a natural class of equivalence
which is clear from context. Here
these equivalences are: the existence of an equivariant homotopy
between two equivariant maps and of an isotopy
between two embeddings.}
Next, $\Lambda f$ induces an equivariant map
$(\Lambda f)_*: \Til{N}\to S^{m-1}$,
$[(\Lambda f)_*]=\alpha\Lambda [f]$.
\footnote{
The cone maps from $\Emb^{2m-1}N_0$
and from the set of individual embeddings
are both denoted by $\Lambda$.
}
The commutativity of the diagram above
is equivalent to the following fact:
the map $(\Lambda f)_*$ is equivariantly homotopic
to the composition
$$\Til{N}\stackrel{\lambda}\to\Sigma\Til{N_0}\stackrel{\Sigma f_*}\to S^{m-1}.$$
The maps $(\Lambda f)_*$ and $(\Sigma f_*)\lambda$ coincide on $\Til {N_0}$,
both send the `vertical' component
$B^n\times N_0$ of $\Til N\setminus\Til {N_0}$ to the upper
hemisphere of $S^{m-1}$
and the `horizontal' component
$N_0\times B^n$ to the lower hemisphere.
Thus $((\Sigma f_*)\lambda) (x)\in S^{m-1}$ and $(\Lambda f)_* (x)\in S^{m-1}$
are not antipodal for each $x\in\Til N$, meaning that
$(\Sigma f_*)\lambda$ and $(\Lambda f)_*$ are equivariantly homotopic.
$\blacksquare$

\smallskip
The map $\alpha$ is one-to-one by the Haefliger-Weber theorem  [Ha63, We67], [Sk06, 5.2 and 5.4].
The map $\alpha_0$ is one-to-one by the
the Haefliger theorem for manifolds with boundary
(see [Ha63, 6.4],
[Sk02, Theorem 1.1$\alpha\bd$]
for the DIFF case and
[Sk02, Theorem 1.3$\alpha\bd$] for the PL case).
\footnote
{The Haefliger-Weber theorem says that the Haefliger-Wu invariant
$\alpha:\Emb^m N\to\pi^{m-1}\Eq(N)$ is one-to-one for $2m\ge 3n+4$
and the Haefliger theorem for manifolds with boundary says that
$\alpha$ is one-to-one if $N$ has $(n-d-1)$-dimensional spine
for $2m\ge 3n+1-d$
in the DIFF category
and
for $2m\ge 3n+2-d$
in the PL category,
$d\ge 0$.
}
Next, $\Sigma$  is one-to-one
by the equivariant version of Freudenthal suspension
theorem [CF60, Theorem 2.5].

It remains to prove that $\lambda^*$
is surjective and
each
preimage ${\lm^*}^{-1}f_0$ 
is
in 1-1 correspondence 
with a quotient set of
$H_{k+1}(N)^\ts$
for each $f_0\in\pi\Eq^{m-1}(\TN)$.
We will need the following
assertion which is
proved below.

\smallskip
{\bf Assertion 1.}
{\it For a $k$-connected $n$-manifold
and the constructed $\lambda$ consider the 
equivariant cohomology groups
$H^*\Eq(\Sg \Til{N_0},\lambda \Til{N})$,
with respect to the
trivial action of $\Z_2$ on $\Z$-coefficients. 
Then

(a) 
for 
$i\le 2k$
we get
$H^{2n-i}\Eq(\Sg \Til{N_0},\lambda \Til{N})=0$;

(b)
$H^{2n-2k-1}\Eq(\Sg \Til{N_0},\lambda \Til{N})\cong H_{k+1}(N)^\ts$.
}

\smallskip
There is a
1-1 correspondence between
$\pi^{m-1}\Eq\TN$ and $\pi^{m-1}\Eq\lm\TN$
since $\lm$ is not injective only on
some cells of dimension $n<m-2$.
We will thus work with
$\pi^{m-1}\Eq\lm\TN$ and $\pi^{m-1}\Eq\TN$ interchangeably.
Take an equivariant map $f_0:\lm\TN\to S^{m-1}$.
It can be extended to an equivariant map
$f_1:\Sg\Tno\to S^{m-1}$
since by Assertion~1(a),
$$H^i\Eq(\Sg\Tno,\lm \TN; \pi_{i-1}S^{m-1})=0\quad\mbox{for each }i.$$
This proves that $\lm^*$ is surjective.

Fix an extension $f_1:\Sg\Tno\to S^{m-1}$
of $f_0$.
Denote by $\pi\Eq^{m-1}(\Sg \Tno,\lm\TN)$
the set  of equivariant extensions
of $f_0$ on $\Sg\Tno$ up to equivariant homotopy
fixed on $\lm\TN$.
Consider the following
diagram.
$$
\begin{CD}
0\\
@AAA\\
\pi\Eq^{m-1}(\Til{N})\\
@A\lambda^* AA\\
\pi\Eq^{m-1}(\Sigma \Til{N_0}) \\
@Aj AA\\
\pi\Eq^{m-1}(\Sg \Tno,\lm\TN)@>d_2>> H^{m-1}\Eq(\Sg \Tno,\lm\TN)\\
\end{CD}
$$
Here

$\bullet$ $j$ is the natural map;

$\bullet$
$d_2$ is the
`degree' map, well-defined and bijective 
by the relative equivariant
version of the 
Hopf--Whitney Theorem
since
$H^r\Eq(\Sigma\Tno,\lm\TN)=0$ for $r\ge 2m$
by Assertion~1(a).

We obtain ${\lm^*}^{-1}f_0=\img j$ is in 1-1 correspondence
with a quotient of $\pi^{m-1}\Eq(\Sg\Tno,\lm\TN)$.
Main Theorem~2(a) now follows from Assertion~1(b).~$\blacksquare$

\smallskip
{\it Proof of Assertion~1(a),(b)}.
Let us introduce some notation.
Take $P\in N\setminus N_0$.
We denote by $\picH$, $\picV$, $\picVD$ the
following submanifolds of $N\times N$, respectively:
$N_0{\times}\{P\}$, $\{P\}{\times} N_0$, $\{P\}{\times} N_0\sqcup \mathrm{diag}\,N_0$,
where $\mathrm{diag}\,N_0$ is the diagonal embedding.
Let $UZ$ denote a regular neighborhood of an embedded $Z\subset \Sigma\Til{N_0}$.
Let $\con_+(\bd U\picH)\subset \Sigma \Til{N_0}$
be the upper cone over $\bd U\picH$, and analogously
denote the lower cone.

We obtain the following chain of isomorphisms for each
$i< n-2$.
\begin{align*}
H^{2n-i}\Eq(\Sg \Tno,\lm\Tno)
\cong&&
\footnotesize{
\mbox{As follows from the definition of }\lambda
}\\
H^{2n-i}\Eq(\Sigma \Til{N_0},\ \Til{N_0}\cup \con_+(\bd U\picH)\cup\con_-(\bd U\picV))
\cong&&
\\
H^{2n-i}\Eq\left(\frac{(\Sigma \Til{N_0})}{\Til{N_0}},\ \frac{\Til{N_0}\cup\con_+(\bd U\picH)\cup\con_-(\bd U\picV)}{\Til{N_0}}\right)
\cong&&
{\footnotesize
\begin{array}{r}
\mbox{By equivariant homeomorphism between }\\
\mbox{the two pairs; the induced involution $g'$ on}\\
(\Sigma \Til{N_0})\vee (\Sigma \Til{N_0}) \mbox{ is described below}\\
\end{array}
}
\\
H^{2n-i}\Eq((\Sigma \Til{N_0})\vee (\Sigma \Til{N_0}),\ \Sigma(\bd U\picH)\vee\Sigma(\bd U\picV))
\cong&&
\footnotesize{
\mbox{Desuspension isomorphism}
}
\\
H^{2n-i-1}\Eq(\Til{N_0}\vee \Til{N_0},\ \bd U\picH\vee \bd U\picV)
\cong&&
{\footnotesize
\begin{array}{r}
\mbox{Because $g'$ is a diffeomorphism from one  }\\
\mbox{component of } (\Sigma \Til{N_0})\vee (\Sigma \Til{N_0}) \mbox{ onto another}
\end{array}}
\\
H^{2n-i-1}(\Til{N_0},\ \bd U\picH)
\cong&&
\footnotesize{
\mbox{Excision}
}\\
H^{2n-i-1}(N^2-U\picVD,\ \picH)
\cong&&
\footnotesize{
\mbox{By exact sequence of pair, } 2n-i-2>n
}\\
H^{2n-i-1}(N^2-U\picVD)
\cong&&
\footnotesize{
\mbox{Poincar\'e duality}
}\\
H_{i+1}(N^2-U\picVD,\bd)
\cong
&&
\footnotesize{
\mbox{Excision}
}\\
H_{i+1}(N^2,\picVD)
\cong&&
\footnotesize{
\begin{array}{r}
\mbox{This pair has homological type of }\\
\mbox{the smash product, see below }
\end{array}}\\
H_{i+1}(N\wedge N)\cong
&&
\footnotesize{
\mbox{By K\"unneth formula}
}\\
\begin{cases}
0,& 0\le i\le 2k\\
(H_{k+1}N)^{\otimes 2},& i=2k+1
\end{cases}
&&.
\end{align*}

The induced involution $g'$ on $(\Sigma \Til{N_0})\vee (\Sigma \Til{N_0})$
is the composition
of the map changing the two components of
$(\Sigma \Til{N_0})\vee (\Sigma \Til{N_0})$
and of the involution $g$ on $\Sigma \Til{N_0}$ applied componentwise.

Here $N\wedge N:=N^2/(N\vee N)$,
where $N\vee N\subset N^2$ is given by the vertical and horizontal
embeddings.
The isomorphism $H_{i+1}(N^2,\picVD)\cong H_{i+1}(N\wedge N)$
is implied by
the following easy fact.
The map
$H_{i+1}N\to H_{i+1}N^2$ induced by diagonal embedding
coincides with the composition
$$H_{i+1}N
\stackrel{\mathrm{id}\oplus \mathrm{id}}\to
H_{i+1}N\oplus H_{i+1}N
\stackrel{\mathrm{v}\oplus \mathrm{h}}\to
H_{i+1}N^2$$
where $\mathrm{v}$, $\mathrm{h}$
are induced by vertical and horizontal embeddings, respectively.
$\blacksquare$

\smallskip
{\it Proof of Conjecture~1 (sketch)}.
We continue from the place where the proof
of Theorem~1 ended. Inside this proof we set $k=0$.
The diagram from above can be completed up
to the following commutative diagram.
$$
\begin{CD}
0\\
@AAA\\
\pi\Eq^{m-1}(\Til{N})\\
@A\lambda^* AA\\
\pi\Eq^{m-1}(\Sigma \Til{N_0})@>d_1>> H^{m-1}\Eq(\Sg \Tno)\\
@Aj AA@A j' AA\\
\pi\Eq^{m-1}(\Sg \Tno,\lm\TN)@>d_2>> H^{m-1}\Eq(\Sg \Tno,\lm\TN)\\
\end{CD}
$$

Here the new map
$d_1$ is the `degree' map
(i.e. the first obstruction for a given map to be
equivariantly homotopic to $f_1$),
well-defined and bijective 
by the equivariant
version of the 
Hopf--Whitney Theorem [Pr06,  p.~103 Theorem 10.5]
since
$H^r\Eq(\Sigma\Tno)=0$ for $r\ge 2m$
[Sk10, Deleted Product Lemma].

The preimage ${\lm^*}^{-1}f_0=\img j$
is in 1-1 correspondence with $\img j'$, and
Conjecture 1 follows from the following Assertion
(we need only the case $k=0$).
$\blacksquare$

\smallskip
{\bf Assertion 1} 
{\it 
(c).
For a $k$-connected $n$-manifold $N$
and the constructed $\lambda$ consider the 
equivariant cohomology groups
$H^*\Eq(\Sg \Til{N_0},\lambda \Til{N})$,
$H^*\Eq(\Sg \Til{N_0})$
with respect to the
trivial action of $\Z_2$ on $\Z$-coefficients. 
Then
the image of
$j':\ H^{2n-2k-1}\Eq(\Sg \Til{N_0},\lambda \Til{N})
\to H^{2n-2k-1}\Eq(\Sg \Til{N_0})$
is isomorphic to $B_{n+k} (H_{k+1}(N))$.
}

\smallskip
{\it Proof of Assertion~1(c) (sketch)}.
We continue from the place where the proof
of Assertion~1(a),(b) ended.
Let $\phi:H_{k+1}(N)^{\otimes 2}\to H^{2n-2k-1}\Eq(\Sg\Tno,\lm\TN)$
be the inverse to the composition the 
long chain of isomorphisms
above.
Then $\phi=l\sigma m$
on the following commutative diagram.
$$
\begin{CD}
H_{k+1}(N)^\ts
\\
@VmV\cong V
\\
H^{2n-2k-2}(\Tno,\bd U\picH)
@>r>\subset >
H^{2n-2k-2}(\Tno)
\\
@V\sigma V \cong \quad x\mapsto \Sigma x\vee \Sigma gxV
@V\bar\sigma V\cong \quad x\mapsto \Sigma x\vee \Sigma gxV
\\
H^{2n-2k-1}\Eq((\Sigma \Til{N_0})\vee (\Sigma \Til{N_0}),\ \Sigma(\bd U\picH)\vee\Sigma(\bd U\picV))
@>>>
H^{2n-2k-1}\Eq((\Sigma \Til{N_0})\vee (\Sigma \Til{N_0}))
\\
@VlV\cong V
@V\bar lV\twoheadrightarrow \quad x\vee y\mapsto x+(-1)^{n+1}yV
\\
H^{2n-2k-1}\Eq(\Sigma \Til{N_0},\lambda\Til{N_0})
@>j'>>
H^{2n-2k-1}\Eq(\Sigma \Til{N_0})
\end{CD}
$$
Here
$m,\sigma,l$
are the isomorphisms from the chain of isomorphisms above,
horizontal maps are natural
$\bar\sigma$ and $\bar l$ are the 
non-relative analogues of $\sigma$ and $l$.
Let $g:N^2\to N^2$ be the standard involution.
By $g$ we will also denote induced maps in 
(co)homology groups
induced by this involution.
Clearly, $\sigma(x)=\Sigma x \vee \Sigma gx$
The analogous formula holds for $\bar\sigma$.
We also get
$\bar l(x\vee y)=x+(-1)^{n+1}y$
since the map $g$ has degree $(-1)^{n+1}$.
\footnote{Note that $\sigma,\bar\sigma, l$ are isomorphisms
$\bar l$ is epimorphic.
We will now prove $r$ is monomorphic
and it will follow that
$\Ker \bar l=\Ker j'$.}

Let us show that $r$ is monomorphic.
Indeed, $\Ker r$ is the image of $H^{2n-2k-1}(\Tno)$
under the map from exact sequence of pair.
But $H^{2n-2k-1}(\Tno)\cong H_{2k+1}(N^2,\picHVD)\cong 0$
where 
$\picHVD=N_0{\times}\{P\}\sqcup \{P\}{\times} N_0\sqcup\mathrm{diag}\,N_0$
and the last isomorphism is analogous
to the isomorphism $H_{2k+1}(N^2,\picVD)\cong 0$
proved above.

Suppose $z\in H_{k+1}(N)^\ts$.
Then $z\in\Ker \phi$
if and only if $\Sigma rm(z+(-1)^{n+1}gz)=0$.
Since $r$ is monomorphic, this is
equivalent to
$z=(-1)^{n}gz$.
Finally,
for $a,b\in H_{k+1}(N)^\ts$
we get
$g(a\otimes b)=(-1)^{k+1}(b\otimes a)$,
so
$z=(-1)^{n}gz$
if and only if
$z\in B_{n+k+1}(H_{k+1}(N))$.
Thus
$\Ker j'=B_{n+k+1}(H_{k+1}(N))$
and so
$\img j'=B_{n+k}(H_{k+1}(N))$.
Part~(b) of Assertion~1 is proved.~$\blacksquare$

\smallskip
{\bf Acknowledgements.}
The author is grateful to 
D.~Crowley,
D.~Gon\c calves,
S.~Melikhov
for useful discussions
and especially to
A.~Skopenkov for constant help and support.

\end{document}